\documentclass[twoside]{article}

\usepackage{amsmath,amsthm,amssymb}

\usepackage{newtxtext,newtxmath}

\usepackage[text={12.5cm,19cm},centering,paperwidth=17cm,paperheight=24cm]{geometry}

\usepackage[fontsize=10.4pt]{scrextend}

\def\titlerunning#1{\gdef\titrun{#1}}

\makeatletter
\def\author#1{\gdef\autrun{\def\and{\unskip, }#1}\gdef\@author{#1}}
\def\address#1{{\def\and{\\\hspace*{15.6pt}}\renewcommand{\thefootnote}{}\footnote{#1}}\markboth{\autrun}{\titrun}}
\makeatother

\def\email#1{email: \href{mailto:#1}{#1} }
\def\subjclass#1{\par\bigskip\noindent\textbf{Mathematics Subject Classification 2020.} #1}
\def\keywords#1{\par\smallskip\noindent\textbf{Keywords.} #1}

\newenvironment{dedication}{\itshape\center}{\par\medskip}
\newenvironment{acknowledgments}{\bigskip\small\noindent\textit{Acknowledgments.}}{\par}

\newtheorem{thm}{Theorem}[section]
\newtheorem{cor}[thm]{Corollary}
\newtheorem{lem}[thm]{Lemma}

\newtheorem{prop}[thm]{Proposition}

\theoremstyle{definition}

\newtheorem{exa}[thm]{Example}

\newtheorem*{rem}{Remark}

\numberwithin{equation}{section}

\usepackage[hyperfootnotes=false,colorlinks=true,allcolors=blue]{hyperref}

  \newcommand{\slim}{\operatorname{s-lim}}

 \newcommand{\supp}{\operatorname{supp}}
  \newcommand{\e}{\eqref}

\newcommand{\q}{\quad}

\newcommand{\ti}{\tilde}
\newcommand{\wt}{\widetilde}

 \renewcommand{\d}{\delta}

  \newcommand{\ran}{\operatorname{Ran}}

\newenvironment{pf}{\begin{proof}}{\end{proof}}

\def\qqq{\mathrel{\subset\mkern-15mu\lower.38ex\hbox{${\scriptscriptstyle\rightarrow}$}}}

\let\cal\mathcal

\let\Bbb\mathbb

\begin{document}

\titlerunning{Scattering theory for Laguerre operators}

\title{\textbf{Scattering theory for Laguerre operators}}

\author{D. R. Yafaev}

\date{}

\maketitle

\address{Univ  Rennes, CNRS, IRMAR-UMR 6625, F-35000
    Rennes, France and SPGU, Univ. Nab. 7/9, Saint Petersburg, 199034 Russia; 
    \email{yafaev@univ-rennes1.fr}}

\begin{dedication}
To Ari Laptev on the occasion of his 70th birthday
\end{dedication}


\begin{abstract}
We study Jacobi operators $J_{p}$, $p> -1$, whose eigenfunctions are Laguerre polynomials. 
All operators $J_{p}$ have absolutely continuous simple spectra coinciding with the positive half-axis. This fact, however, by no means imply that the wave operators for the pairs $J_{p}$, $J_{q}$  where $p\neq q$ exist.  Our goal is to show that, nevertheless, this is true and to find explicit expressions for these wave operators. We also study the time evolution of $(e^{-J  t} f)_{n}$  as $|t|\to\infty$ for Jacobi operators $J$ whose eigenfunctions are different classical polynomials.  For Laguerre polynomials,
it turns out that  the evolution $(e^{-J_{p} t} f)_{n}$ is concentrated in the region where $n\sim t^2$ instead of $n\sim |t |$ as happens in standard situations. 

As a by-product of our considerations, we
obtain  universal relations between amplitudes and phases in asymptotic formulas for general orthogonal polynomials.
   \end{abstract}

\subjclass[2000]{33C45, 39A70, 47A40, 47B39}

\keywords{ Jacobi operators, Laguerre polynomials, asymptotic formulas}  





\maketitle


\section{Introduction}  

\subsection{Jacobi operators} 

Jacobi operators $J$ are discrete analogues of  one-dimensional differential   operators. They are defined in the space $ \ell^2 ({\Bbb Z}_{+})$  by the formula
  \begin{equation}
(J f)_{n}= a_{n-1} f_{n-1}  +b_{n} f_{n }   + a_{n} f_{n+1} ,  \q n\in{\Bbb Z}_{+}, \q a_{-1}=0. 
\label{eq:RR}\end{equation}
We always suppose that $a_{n }>0$ and $b_{n}=\bar{b}_{n}$. Of course spectral properties of 
Jacobi operators depend crucially on a behavior of the coefficients $a_{n}$ and $b_{n}$ as $n\to\infty$. 

 In the simplest case $a_{n}= a_{\infty}>0$, $b_{n}=0$,   the operator $J$ (known as the ``free" discrete Schr\"odinger operator) has the absolutely continuous spectrum $[-2 a_{\infty}, 2a_{\infty}]$, and its eigenfunctions are expressed via Chebyshev polynomials of second kind.  Operators $J$ whose eigenfunctions are Jacobi polynomials are natural generalizations of this operator.

The situation $a_{n}\to\infty$ as $n\to\infty$ is also quite common in applications. We discuss only the case  where the   Carleman condition 
    \begin{equation}
    \sum_{n=0}^\infty a_{n}^{-1}=\infty
\label{eq:Carl}\end{equation}
is satisfied.
Suppose that    $b_{n}/2 a_{n} \to \gamma$  as $n\to\infty$.  If $|\gamma|< 1$, then the spectrum of the operator $J$ is purely absolutely continuous and coincides with the whole axis $\Bbb R$; see \cite{Jan-Nab, Apt}.  A famous example is
 \begin{equation}
      a_{n}=\sqrt{ (n+1)/2} , \q b_{n}=0. 
     \label{eq:H}\end{equation}
     Eigenfunctions of the corresponding operator $J$ are the Hermite polynomials.
If $|\gamma|> 1$ ($| \gamma | =\infty$ is admitted), then the spectrum of $J$ is discrete.  The   case $|\gamma| =1$ is critical, and the spectral properties of $J$  depend on details of the behavior of $b_{n}/(2 a_{n} ) -\gamma$  as $n\to\infty$. The results pertaining to this situation are scarce. We mention only papers \cite{Jan-Nab-Sh, Na-Si} and references therein.

Here we study an important particular case  
 \begin{equation}
    a_{n} =  a_{n}^{(p)}= \sqrt{(n+1)(n+1+p)} \q\mbox{and}\q     b_{n} = b_{n}^{(p)}=  2n+p+1, \q p>-1. 
\label{eq:Lag}\end{equation}
Thus we have   $\gamma=1$.
For all $p$, the Jacobi operators   $J=J_{p}$   (we call them the Laguerre operators)   with the recurrence coefficients \e{eq:Lag}
  have absolutely continuous  spectra
coinciding with $[0,\infty)$. Eigenfunctions of these operators are the Laguerre polynomials $L_{n}^{(p)} (\lambda)$. 
Our goal is to  investigate an asymptotic behavior of the unitary group $e^{-i J_{p} t}$ as $t\to\pm \infty$.  We show that, for different $p$, these asymptotics are essentially the same although the operators $J_{q}  - J_{p}$ are not even   compact unless $q=p$.

Another goal of the paper is to obtain detailed asymptotic formulas for $(e^{-i J  t}f)_{n}$ as $| t |\to  \infty$  for sufficiently arbitrary  Jacobi operators $J$.  Here we suppose that  an
  asymptotic behavior of the corresponding orthogonal polynomials $P_{n}  (\lambda)$  as $n\to  \infty$ is known. These general results are illustrated on examples  of the classical polynomials.    In particular, for Laguerre operators $J_{p}$, we show that    the evolution $(e^{-J_{p} t} f)_{n}$ is concentrated in the region where $n\sim t^2$ instead of $n\sim |t |$ as happens in standard situations.




\subsection{Scattering theory}  

We work in the framework of scattering theory. Let us briefly  recall  its basic notions. We refer, e.g., to the book \cite{Ya}  for a   detailed presentation. Consider a couple of self-adjoint operators $A$ and $B$ in some Hilbert space ${\cal H}$. In view of our applications, we suppose that both of these operators are absolutely continuous.  Scattering theory  studies the strong limits
\begin{equation}
    W_{\pm} =  W_{\pm} (B,A)= \slim_{t\to\pm \infty}  e^{iB t} e^{-iA t}
\label{eq:WO}\end{equation}
known as the wave operators for the pair $A$, $B$. If the limits \e{eq:WO} exist, then the wave operators   possess several useful features. In particular, they are isometric and  enjoy the intertwining property
$
BW_{\pm} =W_{\pm} A .
$
It follows that the restriction of the operator $B$ on the image $\ran W_{\pm}$ of $W_{\pm}$ is unitary equivalent to the operator $A$.  If both wave operators $W_{\pm} (B,A)$ and $W_{\pm} (A,B)$ exist, then the operators $A$
and $B$ are unitarily equivalent. In this case the spectra of the  operators $A$
and $B$ coincide. 

The existence of the limits  \e{eq:WO} is a non-trivial problem. We emphasize that  the unitary equivalence of operators $A$ and $B$ does not imply the existence of the wave operators \e{eq:WO}. A notorious example is given by the pair of multiplication $A  $, 
$(Au)(x)=  x u(x)$, and differential   $B=-id /dx$ operators in the space $L^2 ({\Bbb R})$.  Generally speaking, the limits 
 \e{eq:WO} (or their appropriate modifications) exist if the perturbation $B-A$ is in some sense small. This might mean different things.  For example, it suffices to assume that the operator $B-A$ is trace class or that it acts as an integral operator with smooth kernel in the diagonal representation of the operator $A$ (or $B$).
 
 
 Our aim  is to develop scattering theory for  pairs of the operators $J_{p}$, $J_{q}$.  According to \e{eq:Lag} we have
  \begin{equation}
   a_{n}^{(q)}-a_{n}^{(p)}=  (q-p)/2  + O(n^{-1} ), \q n \to \infty, \q \mbox{and} \q     b_{n}^{(q)}-b_{n}^{(p)}=  q-p, 
\label{eq:Lag1}\end{equation}
so that the operator $J_{q} - J_p$ is not even compact.  Therefore standard methods of scattering theory do not work in this case. 
It turns out however that, from the viewpoint of scattering theory, the diagonal and off-diagonal terms in \e{eq:Lag1} compensate each other.  Note that if only one of the coefficients $a_{n}$ or $b_{n}$ is changed, then the spectrum of the operator  $J_{p}$  is shifted.  In this case the  wave operators cannot exist.

Although the difference $J_{q} - J_p$ is  by no means small, there exists a natural one-to-one correspondence between eigenfunctions of the operators $J_{p}$ for different $p$. Their asymptotics as $n\to\infty$ differ by a phase shift only.  This allows us to show that the wave operators $W_{\pm}  (J_{q}, J_{p})$ exist for all $p, q> -1$. This result is obtained by a direct calculation which yields also     explicit expressions for  the wave operators.

     \subsection{Structure of the paper}
     
     Jacobi operators  and associated orthogonal polynomials, in particular, Laguerre polynomials,  are discussed in Sect.~2.
     
      In Sect.~3, we prove the existence of the wave operators  $W_{\pm}  (J_{q}, J_{p})$ (Theorem~\ref{WW}). We also construct scattering theory for pairs $J$, $\wt{J}$ where $J=J_{p}$ for some $p> -1$ and the coefficients of a Jacobi operator $\wt{J}$ are sufficiently close to those of $J$  (Theorem~\ref{PT}). 
     
     In Sect.~4.1,   we exhibit  a link between a large time behavior of $( e^{-i J  t}f )_{n}$ and asymptotics of associated orthogonal polynomials $P_{n}  (\lambda)$ as $n\to\infty$.   This leads to universal  relations between amplitudes and phases in asymptotic formulas for $P_{n}  (\lambda)$.  These results are illustrated in    Sect.~4.2 -- 4.4 on examples of classical polynomials. The results for Laguerre, Jacobi and Hermite polynomials are stated 
       as Theorems~\ref{WU}, \ref{WJJ} and \ref{WH}, respectively. The Hermite operator  $J$ is somewhat exceptional since the evolution $e^{-i J t}f $  is dispersionless in this case.

     Below $\| \cdot \|$ is the norm in the space $ \ell^2  ({\Bbb Z}_{+})$; $I$ is the identity operator; $C$ are different positive constants whose precise values are of no importance.
     

\section{Jacobi operators  and orthogonal polynomials}  

    
      \subsection{Orthogonal polynomials}
    
    Here we collect  necessary information about the Jacobi operators $J$ given by formula \e{eq:RR} and associated orthogonal polynomials $P_{n } (z)$; see    the books \cite{AKH, Nev}
 for a comprehensive presentation.
        Given coefficients $a_{n}> 0$, $b_{n}= \bar{b}_{n}$, $n\in{\Bbb Z}_{+}$, one constructs
  $P_{n } (z)$   by the recurrence relation
 \begin{equation}
 a_{n-1} P_{n-1} (z) +b_{n} P_{n } (z) + a_{n} P_{n+1} (z)= z P_n (z), \q n\in{\Bbb Z}_{+}, 
 \q z\in{\Bbb C}, 
\label{eq:Py}\end{equation}
and the boundary conditions $P_{-1 } (z) =0$, $P_0 (z) =1$. Then $P_{n } (z)$ is a polynomial of degree $n$.
Obviously, $P(z)=\{ P_{n} (z)\}_{n=0}^\infty$ satisfies the equation $J P(z)=zP(z)$, that is, it is an ``eigenvector"  of the operator $J$.

 
   We consider the operator $J$  in the space $\ell^2 ({\Bbb Z}_{+})$. Let us denote by $J_{0}$ 
the minimal   operator    defined    by  formula \e{eq:RR} on a set $\cal D  $ of vectors $f=\{ f_{n}  \}_{n=0}^\infty$ with only a  finite number of non-zero components.  This operator is symmetric; moreover, it is essentially self-adjoint   if the   Carleman condition 
 \e{eq:Carl} 
is satisfied. In particular,   condition   \e{eq:Carl}  holds true for the coefficients  \e{eq:Lag}. For essentially self-adjoint operators, the domain ${\cal D} ( J)$ of the  closure $J$ of the operator $J_{0}$ consists of all vectors  $f \in \ell^2 ({\Bbb Z}_{+})$ such that $J f \in \ell^2 ({\Bbb Z}_{+})$.
   The spectrum of  the self-adjoint operator $J$ is  simple with $e_{0} = (1,0,0,\ldots)^{\top}$ being a generating vector. Therefore it is natural to define   the   spectral measure of $J$ by the relation $d\rho (\lambda)=d(E (\lambda)e_{0}, e_{0})$ where  $E (\lambda)$      is the spectral family of the operator $J$. The  polynomials  $P_{n}(\lambda)$ (we call them orthonormal) are orthogonal and normalized  in the spaces $L^2 ({\Bbb R};d\rho )$: 
    \begin{equation}
\int_{-\infty}^\infty P_{n}(\lambda) P_{m}(\lambda) d\rho (\lambda) =\d_{n,m};
\label{eq:J5}\end{equation}
as usual, $\d_{n,n}=1$ and $\d_{n,m}=0$ for $n\neq m$.

Alternatively,   given a probability measure $d\rho(\lambda)$, the polynomials $P_0 (\lambda),P_1 (\lambda),\ldots, \\ P_{n} (\lambda),\ldots$ can be obtained by the Gram-Schmidt orthonormalization of the monomials $1,\lambda,\ldots,\lambda^n,\ldots$ in the space $L^2({\Bbb R}_{+}; d\rho)$. It is easy to see that $P_{n} (\lambda)$ is a polynomial of degree $n$, that is, 
$P_{n} (\lambda)=k_{n}(\lambda^n+ r_{n}\lambda^{n-1} +\cdots)$ with $k_{n} \neq 0$. 
One usually  requires     $k_{n} >0$.   The recurrence coefficients $a_{n}, b_{n}$ can be recovered by the formulas  $a_{n}= k_{n}  / k_{n+1}$, $b_{n}= r_{n}
 - r_{n+1}$.
 
 One  defines a mapping $\boldsymbol{\Phi}: \ell^2 ({\Bbb Z}_{+})\to L^2 ({\Bbb R}; d\rho)$ by the formula 
    \begin{equation}
 (\boldsymbol\Phi f)(\lambda)=\sum_{n=0}^\infty P_{n} (\lambda) f_{n}, \q f=\{ f_{n}\}_{n=0}^\infty \in {\cal D}.
\label{eq:UU}\end{equation}
This mapping is isometric according to \e{eq:J5}. It is also unitary if  the set of all polynomials  $P_n(\lambda) $, $n\in {\Bbb Z}_{+}$,  is dense in $L^2 ({\Bbb R}; d\rho)$.  This condition is satisfied if the operator $J_{0}$ is
 essentially self-adjoint. 
  Putting together definitions \e{eq:RR}, \e{eq:Py}  and \e{eq:UU}, it is easy to check   the intertwining property 
$(\boldsymbol{\Phi} J f) (\lambda)= \lambda (\boldsymbol{\Phi} f)(\lambda)$.
 
  Typically, on the absolutely continuous spectrum of a Jacobi operator $J$  when $d\rho(\lambda)= \tau(\lambda) d\lambda$,
  the orthonormal polynomials $P_{n} (\lambda)$  have oscillating asymptotics   
     \begin{equation} 
       \boxed{
P_{n}  (\lambda) = 2 \kappa   (\lambda)n^{-r} \cos \Omega_{n}(\lambda)  + o(n^{-r}), }
\label{eq:HX1}\end{equation}
  \begin{equation} 
    \boxed{
  \Omega_{n}'(\lambda) = \omega(\lambda)n^s + o(n^s) }
\label{eq:HX2}\end{equation}
where $ \kappa   (\lambda)> 0$, $s>0$ and we can suppose that  $ \omega  (\lambda)> 0$.
The exponents $r$,  $s$,  the amplitude  $ \kappa   (\lambda)$ and the phase $\omega(\lambda)$ are determined by the recurrence coefficients $a_{n}$ and $b_{n}$. Our considerations (see Sect.~4.1) show that these quantities are necessarily  linked by universal  relations:
  \begin{equation} 
  \boxed{
2 r+ s=1,}
\label{eq:HXG}\end{equation}
\begin{equation} 
  \boxed{
2 \pi \tau(\lambda) \kappa^2 (\lambda) = s   \omega (\lambda) .}
\label{eq:HXG1}\end{equation}

        \subsection{Laguerre operators}
 
   Suppose now that  the recurrence coefficients $a_{n}$, $b_{n}$ are given  by formulas \e{eq:Lag}. In this case the orthogonal polynomials $L_{n}^{(p)} (z)$ defined by     relations  \e{eq:Py} with the   boundary conditions $L^{(p)}_{-1 } (z) =0$, $L^{(p)}_{0 } (z) =1$ are known as the Laguerre polynomials.  
  Note that the normalized polynomials $L_{n}^{(p)} (z)$ we consider here are related to the Laguerre polynomials ${\bf L}_{n}^{(p)} (z)$ defined in \S 10.12 of the book \cite{BE} or in \S 5.1 of the book \cite{Sz} by the equality
\begin{equation}
L_{n}^{(p)} (z)= (-1)^n   \sqrt{\frac{\Gamma (1+n) \Gamma (1+p)}
{\Gamma (1+n+p)  }} \:{\bf L}_{n}^{(p)} (z).
\label{eq:Lag3}\end{equation}
According to asymptotic formula (10.15.1) in \cite{BE} for positive $\lambda$, we have
\begin{equation}
L_{n}^{(p)} (\lambda)=
(-1)^n   \sqrt{\frac{\Gamma (1+p)  }
 {\pi  }}  \lambda^{-p/2-1/4}  e^{\lambda/2} n^{-1/4} \cos \Big( 2\sqrt{n\lambda}-\frac{2p+1} {4} \pi\Big) +O(n^{-3/4})  
\label{eq:Lag4}\end{equation}
as $n\to\infty$.
This asymptotics is uniform in $\lambda\in [\lambda_{0}, \lambda_{1}]$  if $0< \lambda_{0} <\lambda_{1}< \infty$.
 
            Let us now consider the Laguerre operators $ J_{p}$ defined by formula \e{eq:RR} where $a_{n}$, $b_{n}$ are  given by \e{eq:Lag}. The spectral measures of the operators $ J_{p}$ are supported on the half-axis  $[0,\infty)$, they are absolutely continuous and are given by the relation (see, e.g., formula (10.12.1) in \cite{BE}) 
         \begin{equation}
d\rho_{p} (\lambda)=  \tau_{p}(\lambda)d\lambda  \q \mbox{where}\q \tau_{p}(\lambda)=\frac{1}{\Gamma(p+1)}\lambda^p e^{-\lambda}, \q \lambda\in {\Bbb R}_{+}.
\label{eq:Lag2}\end{equation}

Since the  measure $d\rho_{p} (\lambda)$  is absolutely continuous, it is convenient to reduce the Jacobi operator $J_{p}$ to the operator  ${\sf A}$ of multiplication by $\lambda$ in the space $L^2 ({\Bbb R}_+)$.  To that end, we put
   \begin{equation}
\varphi_{n}^{(p)}(\lambda)= \sqrt{\tau_{p}(\lambda)}  \: L_{n}^{(p)} (\lambda), \q \lambda\in {\Bbb R}_{+},
\label{eq:UF}\end{equation}
  and introduce a mapping $\Phi_{p} : \ell ^2 ({\Bbb Z}_{+})\to L^2 ({\Bbb R}_{+} )$  by the formula (cf. \e{eq:UU})
  \begin{equation}
 (\Phi_{p} f)(\lambda)=\sum_{n=0}^\infty\varphi_{n}^{(p)}  (\lambda) f_{n}, \q f=\{ f_{n}\}_{n=0}^\infty \in {\cal D}, \q \lambda\in {\Bbb R}_{+}.
\label{eq:UU1}\end{equation}
The operator $\Phi_{p}^* :  L^2 ( {\Bbb R}_{+} )\to  \ell^2 ({\Bbb Z}_{+})$ adjoint to $\Phi_{p} $  is given by the equality
   \[
(\Phi_{p}^* g)_{n}=   \int_0^\infty \varphi_{n}^{(p)} (\lambda ) g(\lambda) d\lambda,\q n\in {\Bbb Z}_{+}.
\]
The operator $\Phi_{p}$ is unitary, that is,
   \begin{equation}
\Phi_{p}^*  \Phi_{p}=I,\q  \Phi_{p}  \Phi_{p}^*   =  I ,
 \label{eq:ps3}\end{equation}
and enjoys  the intertwining property  
  \begin{equation}
 \Phi_{p} \, J_{p} = {\sf A} \Phi_{p}  .
 \label{eq:inw1}\end{equation} 

 
 \section{Wave operators}
 
  \subsection{Two Laguerre operators}
 
 One of our main results is stated as follows.
 
  \begin{thm}\label{WW}
   Let $J_{p}$ be  a Jacobi operator with matrix elements \e{eq:Lag} in the space $ \ell^2  ({\Bbb Z}_{+}) $.  Define the unitary  operators $\Phi_p: \ell^2  ({\Bbb Z}_{+}) \to L^2  ({\Bbb R}_{+}) $   by formulas \e{eq:UF} and \e{eq:UU1}.  Then for all $p,q>-1$, the wave operators $W_{\pm} (J_{q} , J_{p})$ exist and 
     \[
W_{\pm} (J_{q} , J_{p})= e^{\pm i (q -p)\pi/2} \Phi_{q}^{*}  \Phi_{p} .
\]
 \end{thm}
 
 We start a proof with a simple standard statement.   
 
   \begin{lem}\label{WW1}
   The claim of Theorem~\ref{WW} is equivalent to the relation
   \begin{equation}
\slim_{t\to\pm \infty} (\Phi_p^{*}-   \mu_{\pm} \Phi_q^{*} ) e^{-i{\sf A} t}= 0, \q \mu_{\pm}= e^{\pm i (q-p)\pi/2} .
\label{eq:WW1}\end{equation}
    \end{lem}
    
      \begin{pf}
      Let $f$ be an arbitrary element of the space $\ell ^2 ({\Bbb Z}_{+})$ and $g = \Phi_{p}f$.
      In view of the properties  \e{eq:ps3} and \e{eq:inw1}, we have
\[
\| e^{iJ_{q}t}e^{-iJ_{p}t} f - \mu_{\pm} \Phi_{q}^{*}  \Phi_{p}f \| = 
\|  e^{-iJ_{p}t} \Phi_{p}^{*} g  - \mu_{\pm} e^{-iJ_q t} \Phi_{q}^{*}   g \| = 
\|  (\Phi_{p}^{*}  - \mu_{\pm} \Phi_{q}^{*} )e^{-i{\sf A}t} g \|. 
\]
Since  the left- and right-hand sides here tend to zero at the same time, this concludes the proof.
    \end{pf}

It suffices to  check  \e{eq:WW1} on a set $C_{0}^\infty ({\Bbb R}_{+})$ dense in $L^2({\Bbb R}_{+})$ .    Let the function  $\varphi_{n}^{(p)} (\lambda )  $ be defined by equalities  \e{eq:Lag2} and \e{eq:UF}. It follows from asymptotic formula \e{eq:Lag4}   that
    \begin{equation}
\varphi_{n}^{(p)} (\lambda )  =
(-1)^n  2^{-1} \pi^{-1/2}  \lambda
   ^{-1/4}   (    n+1)^{-1/4} \big( \nu_{p} e^{2i\sqrt{n\lambda} } + \bar\nu_{p} e^{-2i\sqrt{n\lambda} }\big)    +  r_{n}^{(p)} (\lambda)
   \label{eq:WW2}\end{equation}
where  $\nu_{p}=e^{- i \pi (2p+1)/4 }$ and 
 \begin{equation}
   | r_{n}^{(p)} (\lambda)| \leq C ( n+1)^{-3/4}  
\label{eq:WW3}\end{equation}
uniformly on every compact subinterval of ${\Bbb R}_{+}$.  

 Let us define mappings $V_{\pm} :
C_{0}^\infty ({\Bbb R}_{+}) \to \ell^2 ({\Bbb Z}_{+})$  by the formula
 \begin{equation}
( V_{\pm} g)_{n}  = (-1) ^n 2^{-1} \pi^{-1/2}  (n+1)^{-1/4}
 \int_{0}^\infty e^{\pm  2 i \sqrt{n\lambda} } \lambda^{-1/4} g(\lambda) d\lambda.
   \label{eq:WW4}\end{equation}
  Equality \e{eq:WW2}  implies that
    \begin{equation}
    \Phi_{p}^* e^{-i{\sf A}t}g =   \nu_{p}  V_{+} e^{-i{\sf A}t}g
    + \bar\nu_{p}  V_{-} e^{-i{\sf A}t}g  + R_{p}  (t) g
   \label{eq:WW6}\end{equation}
   where  
   \[
   (R_{p}  (t) g)_{n} =   \int_{0}^\infty  r_{n}^{(p)} (\lambda)  e^{-i\lambda t}
  g(\lambda) d\lambda.
   \]
   
   First, we check that the remainder in \e{eq:WW6} is negligible.
   
    \begin{lem}\label{WW2a}
     Let $g\in C_{0}^\infty ({\Bbb R}_{+})$. Then 
     \begin{equation}
    \lim_{| t|\to \infty} \| R_{p}  (t) g\| =0.
   \label{eq:WW61}\end{equation}   
    \end{lem}
    
     \begin{pf}
     By the Riemann-Lebesgue lemma, every integral in the sum 
 \begin{equation}
\| R_{p}  (t) g \|^2= \sum_{n=0}^\infty \Big| \int_{0}^\infty  r_{n}^{(p)}  (\lambda) e^{- i\lambda t} g(\lambda) d\lambda\Big| ^2
\label{eq:WW5}\end{equation}
tends to zero  as $ |t| \to\infty$.  Estimate \e{eq:WW3} allows us to use the dominated convergence theorem.  Therefore the sum \e{eq:WW5} tends to zero as $ |t| \to\infty$. 
        \end{pf}


 Note a formula  
 \begin{equation}
       \int_{0}^\infty e^{ \pm 2 i \sqrt{n\lambda} -i\lambda t} G(\lambda) d\lambda
       =     i \int_{0}^\infty e^{\pm  2 i \sqrt{n\lambda} -i\lambda t} \Big( \frac{G(\lambda)}  {\pm \sqrt{n}\lambda^{-1/2}-t} \Big)'d\lambda
    \label{eq:WW7A}\end{equation}
        which       can be  verified  by a direct integration by parts.  Using it, we obtain  the following elementary assertion.

     \begin{lem}\label{WW2}
     For an arbitrary $G\in C_{0}^\infty ({\Bbb R}_{+})$,  an estimate
   \begin{equation}
\Big|\int_{0}^\infty e^{ \pm 2 i\sqrt{n\lambda}- i\lambda t} G(\lambda) d\lambda \Big|\leq C_{k} (\sqrt{n}+ | t| )^{-k}, \q \mp t>0,
\label{eq:WW7}\end{equation}
is true for all $k\in {\Bbb Z}_{+}$.  
    \end{lem}
    
     \begin{pf}
     Suppose that $\supp G\subset [\lambda_{1}, \lambda_{2}]$.
    Then
       \[
        \sqrt{n}\lambda^{-1/2}+ | t | \geq \sqrt{n} \lambda_{2}^{-1/2} + | t |
       \]
       and $\lambda^{-3/2}\leq \lambda^{-3/2}_{1}$.  Therefore the right-hand side of \e{eq:WW7A} is estimated by $ C_1(\sqrt{n}+ | t| )^{-1}$
    which proves   \e{eq:WW7} for $k=1$. Further
       integrations  by parts   in \e{eq:WW7A}, yield  \e{eq:WW7} for an arbitrary $k$.
        \end{pf}

      \begin{cor}\label{WW2+}
     For the operators  \e{eq:WW4} and all $g\in C_{0}^\infty ({\Bbb R}_{+})$, we have
  \begin{equation}
   \lim_{t\to\pm \infty} \| V_{\mp}e^{-i{\sf A}t}g \|=0.
\label{eq:WW7+}\end{equation}
    \end{cor}

    Using now relation \e{eq:WW6}  and taking into account Lemma~\ref{WW2a}, we arrive at the following result. 
    
     \begin{lem}\label{WW3}
     Let $g\in C_{0}^\infty ({\Bbb R}_{+})$. Then  
   \begin{equation}
   \lim_{t\to\pm \infty} \|   \Phi_{p}^* e^{-i{\sf A}t}g - \nu_{p}^{\pm 1}V_{\pm}e^{-i{\sf A}t}g \|=0.
\label{eq:WW8}\end{equation}
    \end{lem}
    
    The same result is of course true for $  \Phi_q^* e^{-i{\sf A}t}g$. This  yields relation
   \e{eq:WW1}  with $\mu_{\pm}  = (\nu_{p} \bar{\nu}_{q})^{\pm 1}$. Using Lemma~\ref{WW1}, we conclude the proof of Theorem~\ref{WW}.     $ \q \q\q\q\q\q\Box$
   

   For the scattering operator $S=W_{+}^* W_{-}$, we obtain a very simple expression.
   
    \begin{prop}\label{WWS}
Under the assumptions of Theorem~\ref{WW}, we have
$S =e^{i(p-q)\pi}I$.
    \end{prop}
    
    Recall that the wave operators $W_{\pm} (B,A)$ for  a couple of self-adjoint operators $A$ and $B$  can be represented as  products of an  appropriate Fourier transform for the operator $A$ and the inverse transform  corresponding to $B$. For Schr\"odinger operators in the space $L^2 ({\Bbb R}_{+} )$, this is discussed, for example,  in Section~4.2 (see formula (2.30)) of the book \cite{YA}.  Normally, there are two natural sets of eigenfunctions  of the operators $A$ and $B$.  This leads to two wave operators. In our case these sets of eigenfunctions almost coincide so that the wave operators $W_{\pm} (J_{q}, J_{p})$ differ by a phase factor only whence  the scattering operator is almost trivial.

        \subsection{Perturbation theory}

        Here we choose some $p>-1$ and construct scattering theory for the pair $J=J_{p}$, $\wt{J}=J +V$ where the operator $V$ is in some sense small. We do not assume that $V$ is a Jacobi operator, but, in particular, our results apply to Jacobi operators. We denote by $\wt{\cal H}_{ac}$ the absolutely continuous subspace of the operator $\wt{J}$.
        
         Let us define an operator $N$ in the space  $\ell^2 ({\Bbb Z}_{+})$ by the formula
        \[
        (N f)_{n}=(n+1) f_{n}.
        \]
        Our goal is to prove the following result.
        
        \begin{thm}\label{PT}
   Let $J=J_{p}$, $p> -1$,  be  the Laguerre operator with matrix elements \e{eq:Lag},  and let  $\wt{J}=J +V$ where  $V$ is a symmetric operator  such that  
 \begin{equation}
   V= N^{-r} T N^{-r_{0} } 
   \label{eq:PT}\end{equation}
   for some bounded operator $T$. 
   
   $1^0$ If $  r_{0}>1/4$, $ r>1/4$,
then the wave operators $W_{\pm} (\wt{J}, J)$ exist and are complete, that is,
   $
   \ran W_{\pm} (\wt{J}, J)= \wt{\cal H}_{ac}$.

  $2^0$ If $  r_{0}>1/4$, $ r>1/2$, then the singular spectrum of the operator $\wt{J}$ consists   of eigenvalues of finite multiplicities that may accumulate to the point $0$ only.
    \end{thm}
    
    \begin{cor}\label{PTJ}
   Let $a_{n}$, $b_{n}$,  be  defined by formulas  \e{eq:Lag},  and let  $\wt{J}$ be a self-adjoint Jacobi operator with matrix elements $\ti{a}_{n}$, $\ti{b}_{n}$  such that
   \begin{equation}
   \ti{a}_{n}- a_{n}=O (n^{-\rho}), \q    \ti{b}_{n}-b_{n}=O (n^{-\rho}).
 \label{eq:PT6}\end{equation}
 
 $1^0$   If $\rho>1/2$, then 
   the wave operators $W_{\pm} (\wt{J}, J)$ exist and are complete.
 
  $2^0$  If $\rho>3/4$, then the singular spectrum of the operator $\wt{J}$ consists of eigenvalues   that may accumulate to the point $0$ only.
    \end{cor}
    
 Let us deduce Corollary~\ref{PTJ}  from Theorem~\ref{PT}. We introduce diagonal  matrices $A$ and $B$ with the elements $a_{n}$ and $b_{n}$ and the shift ${\sf S}$: 
 $
 ({\sf S} f)_{n} =f_{n+1}$, $ n\in {\Bbb Z}_{+}$.
   Then 
 $
  J=   A {\sf S}+ {\sf S}^* A + B $
  and with obvious notation, we have
   \[
  \wt{J}-J=   (\wt{A}-A){\sf S}+{\sf S}^* (\wt{A}-A)+ (\wt{B}-B)  .
  \]
  The operators $N^r   (\wt{A}-A) N^{r_{0}}$ and $N^r (\wt{B}-B)  N^{r_{0}}$  are bounded if $r_{0} + r=\rho$. Since the operator $N^{-r} {\sf S} N^{r}$ is also bounded, we see that $N^r ( \wt{J}-J) N^{r_{0}}$ is bounded as long as $r_{0}  + r=\rho$. For the proof of the first statement of Corollary~\ref{PTJ}, we set $r_{0}= r= \rho/2$. The second  statement follows if  $ r_{0}= (\rho-1/4)/2$ and $ r = (\rho+ 1/4)/2$. $\q\q\q \q\q \Box$

    \begin{rem}\label{PTJ2}
    Under the assumptions of  Corollary~\ref{PTJ} one can find asymptotics of the associated orthogonal polynomials. Essentially, it is the same as for the Laguerre polynomials, that is, given (up to a phase shift) by formula \e{eq:Lag4}. 
    \end{rem}
    
      \begin{exa}\label{PTJ3}
  Consider the Jacobi operator $\wt{J}$ with the coefficients
  \[
  \ti{a}_{n}= n+\alpha, \q \ti{b}_{n}= 2n + 2 \alpha-1, \q \alpha >1/2.
  \]
  Up to a shift by $ 2\alpha-1$, this operator is related to the birth and  death processes (see \S 5.2 of the book \cite{Ism}).  Now conditions  \e{eq:PT6} with  $\rho=1$ are satisfied for $a_{n} = a_{n}^{(p)}$, 
  $b_{n} = b_{n}^{(p)}$ where $p= 2(\alpha-1)$. 
    \end{exa}
  
  Note that under the assumptions of Theorem~\ref{PT} or Corollary~\ref{PTJ},  the operator $V=\wt{J}-J  $  belongs to the Hilbert-Schmidt but not to the trace class.  Therefore the assertion about the wave operators $W_{\pm} (\wt{J}, J)$  does not follow from the classical Kato-Rosenblum theorem.
    
        \subsection{Strong  smoothness} 
        
     
     Our proof of Theorem~\ref{PT}
     relies on a notion of strong  smoothness (see Definition~4.4.5 in  the book \cite{Ya}).  We will check strong $J$-smoothness      of the operator $N^{-r}$  for $r>1/4$. 
          Recall that the operator $\Phi$ is defined by formula  \e{eq:UU1} where the functions $\varphi_{n}  (\lambda)$
          are linked to the Laguerre polynomials  $L_{n} (\lambda)$  by equalities  \e{eq:Lag2}, \e{eq:UF}.       
       
     
    \begin{lem}\label{PTJ1}
    Let $\Lambda$ be a compact subinterval of ${\Bbb R}_{+}$ and $r>1/4$.    Then 
               \begin{equation}
| (\Phi N^{-r} f)(\lambda)| \leq C \| f \| .
\label{eq:PT1}\end{equation}
Moreover, 
  \begin{equation}
|(\Phi N^{-r} f)(\mu) - (\Phi N^{-r} f)(\lambda)| \leq C | \mu-\lambda|^s \| f \|
\label{eq:PT2}\end{equation}
 if  $s< 2r -1/2$ and $s\leq 1$.
The constants $C$ in \e{eq:PT1} and \e{eq:PT2} do not depend   on $f\in \ell^2 ({\Bbb Z}_{+})$  and  $\lambda,\mu\in \Lambda$. 
\end{lem}

 \begin{pf}  
It follows from asymptotics  \e{eq:Lag4} of  $L_{n}^{(p)} (\lambda)$ that
  \begin{equation}
|\varphi_{n} (\lambda)| \leq C (1+n)^{-1/4}
\label{eq:PT3}\end{equation}
where the constant $C$ does not depend on $n \in {\Bbb Z}_{+}$ and on $\lambda\in \Lambda$.  So, by the Schwarz inequality, we have
\[
| (\Phi N^{-r} f)(\lambda)|^2 \leq C  \Big( \sum_{n=0}^\infty  (1+n)^{-1/4-r} |f_{n}|\Big)^2
\leq C \sum_{n=0}^\infty  (1+n)^{-1/2- 2r} 
 \| f \|^2 
 \]
 where the series is convergent if $2r>1/2$. This proves \e{eq:PT1}.
 
 For  a proof of \e{eq:PT2}, we need an estimate on derivatives of  $d L_{n}^{(p)} (\lambda)/d\lambda$ of Laguerre polynomials for large $n$.  Let us use  formula (5.1.14) of the book \cite{Sz}  for Laguerre polynomials ${\bf L}_{n}^{(p)} (\lambda)$ linked to $L_{n}^{(p)} (\lambda)$ by equality \e{eq:Lag3}: 
 \[
 \frac{d}{d \lambda} {\bf L}_{n}^{(p)} (\lambda)= -  {\bf L}_{n-1}^{(p+1)} (\lambda)
\q \mbox{whence }
\q
 \frac{d}{d \lambda} L_{n}^{(p)} (\lambda)= \sqrt{\frac{n}{p+1}} L_{n-1}^{(p+1)} (\lambda).
 \]
 It now follows from estimate \e{eq:PT3} that
    \begin{equation}
|\varphi_{n}' (\lambda)| \leq C (1+n)^{1/4} .
\label{eq:PT4}\end{equation}

Note that
\[
|\varphi_{n} (\mu)-\varphi_{n} (\lambda)|  \leq  (2  \sup_{x\in \Lambda}|\varphi_{n} (x)|)^{1-s} \sup_{x\in \Lambda}|\varphi_{n}' (x)|^s  |\mu-\lambda|^s   
\]
for any $s\in [0,1]$.
Therefore using \e{eq:PT3}  and \e{eq:PT4}, we see that
\[
|\varphi_{n} (\mu)-\varphi_{n} (\lambda)|  \leq     C (1+n)^{-1/4 + s/2} |\mu-\lambda|^s .
\]
This yields an estimate
  \begin{multline*}
| (\Phi N^{-r} f)(\mu)-(\Phi N^{-r} f)(\lambda)| \leq \sum_{n=0}^\infty | \varphi_{n}(\mu)-\varphi_{n}(\lambda) |(1+n)^{-r}|f_{n}| 
\\
\leq \| f \|  \sqrt{\sum_{n=0}^\infty | \varphi_{n}(\mu)-\varphi_{n}(\lambda)|^2 (1+n)^{-2r}  }   \leq C |\mu-\lambda|^s \| f \|
\end{multline*}
provided $s< 2r -1/2$.  Thus we get \e{eq:PT2}. 
 \end{pf} 

  The operator   $N^{-r}$ satisfying estimates \e{eq:PT1} and \e{eq:PT2} is called strongly     $J$-smooth  with    exponent $s\in (0,1]$ on the interval      $\Lambda $.

Theorem~4.6.4 of \cite{Ya} states that if a perturbation $V$ admits  representation  \e{eq:PT} with the operators $N^{-r_{0}}$ and $N^{-r }$  strongly     $J$-smooth  (with    some exponents $s_{0} ,s >0$) on all compact subintervals $\Lambda$  of ${\Bbb R}_{+}$, then 
   the wave operators $W_{\pm} (\wt{J}, J)$ exist and are complete.  This is   part~$1^0$ of Theorem~\ref{PT}.
   
   Theorems~4.7.9 and 4.7.10 of \cite{Ya} state that all   spectral results enumerated in part~$2^0$ of Theorem~\ref{PT} are true  provided $s >1/2$ (and $s_{0} >0$).  According to Lemma~\ref{PTJ1} we can choose $s>1/2$ if $r>1/2$. This concludes the proof  of    part~$2^0$ of Theorem~\ref{PT}. $\Box$
   
   Finally, we note that an unusually weak assumption $\rho>1/2$ (instead of the standard $\rho>1$) in \e{eq:PT6})  is explained by a decay \e{eq:PT3} of eigenfunctions of $J$.
   
    \section{Time evolution} 
    
        
    
     It is a common wisdom  that an asymptotic behavior  of $e^{-iJ t} f$ as $t\to\infty$ is determined by spectral properties of the Jacobi operator $J$ and  by asymptotics of the corresponding orthonormal polynomials $P_{n} (\lambda)$ as $n\to\infty$.  
     We first discuss this general idea at a heuristic level and derive new relations between amplitudes and phases in asymptotic formulas for $P_{n} (\lambda)$. Then we illustrate the formulas obtained  on   examples of the classical polynomials.

               \subsection{Universal asymptotic relations} 
     
Assume that the spectrum of a Jacobi operator $J$ is absolutely continuous on an interval $\Lambda$ and the corresponding measure $d\rho(\lambda)=\tau(\lambda)d\lambda$ has a smooth weight $\tau(\lambda)$ for $\lambda\in\Lambda$.  Let
the operator $\boldsymbol{\Phi}$ diagonalizing $J$ be defined by formula \e{eq:UU}. Choose$f$ such that  $f= E(\Lambda)f$  and set
$g(\lambda)=  \sqrt{\tau(\lambda)} (\boldsymbol{\Phi} f)(\lambda)$.  Clearly, $\| g \| _{L^2(\Lambda)}=\| f\|$.
   If  $P_{n} (\lambda)$ satisfy asymptotic relation \e{eq:HX1}, then
    \begin{equation} 
(e^{-i J t} f)_{n}= (n+1)^{-r} \int_{\Lambda} \varkappa   (\lambda)\big(e^{i\Omega_{n}  (\lambda) -i\lambda t}+ e^{-i \Omega_{n}  (\lambda)   -i\lambda t}\big)
 g(\lambda) d\lambda
\label{eq:HX3}\end{equation}
where  $ \varkappa   (\lambda)= \sqrt{\tau(\lambda)} \kappa   (\lambda)$.
   Here and below we keep track only of leading terms in asymptotic formulas. 
  We suppose that the phase $\Omega_{n}  (\lambda)$ obeys condition \e{eq:HX2} where
   $\omega' (\lambda)\neq 0$ for $\lambda\in \Lambda$  and that $g\in C_{0}^\infty (\Lambda)$.

Stationary points of the integrals \e{eq:HX3} are determined by the equations
    \begin{equation} 
\pm \omega (\lambda) n^ s=t.
\label{eq:HX4}\end{equation}
Suppose, for definiteness, that   $t\to+\infty$. Then equation \e{eq:HX4} may have a solution (necessary unique) for the sign $``+"$ only. Let $\sigma= s^{-1}$, $\xi= n/ t^\sigma$, and let $\lambda=\lambda (\xi) $ be the solution of the equation 
   \begin{equation} 
   \omega (\lambda)= \xi^{-s}. 
   \label{eq:HX4A}\end{equation}
  Applying  the stationary phase method to integrals  \e{eq:HX3}, we  see that
    \begin{equation} 
(e^{-i J t} f)_{n}= (2\pi)^{1/2} n^{-r-s/2}   e^{i \phi (\xi,t)}
  |\omega' (\lambda (\xi) )|^{-1/2}
\varkappa   (\lambda (\xi)) g(\lambda (\xi))  
\label{eq:HX5}\end{equation}
where
\[
\phi (\xi,t)=  \pm \pi/4 +\Omega_{\xi t^{\sigma}}(\lambda (\xi))   t - \lambda  (\xi)  t \q\mbox{if}\q \pm\omega' (\lambda)>0.
\]

Let  
 \begin{equation} 
h(\xi)= \xi^{-r-s/2}    |\omega' (\lambda (\xi) )|^{-1/2}
\varkappa   (\lambda (\xi)) g(\lambda (\xi)) 
\label{eq:HX5x}\end{equation}
so that  \e{eq:HX5} reads as
 \begin{equation} 
(e^{-i J t} f)_{n}= (2\pi)^{1/2} t^{-(2r+s)\sigma/2} e^{i \phi (n/t^\sigma,t)}  h(n/t^\sigma).
\label{eq:HX5X}\end{equation}
  Since the operators $e^{-i J t}$  are unitary, it follows from \e{eq:HX5X}  that
 \begin{equation} 
2\pi  \lim_{t\to\infty}  \Big( t^{-(2r+s)\sigma} \sum_{n=0}^\infty   |h ( n/ t^\sigma)|^2 \Big)=\| f\|^2 .
\label{eq:HX5z}\end{equation}
 
  Observe that the integral sums
 \begin{equation} 
  N^{-1} \sum_{n=0}^\infty   |h ( n/ N)|^2\to \int_{\lambda^{-1} (\Lambda)} |h ( \xi)|^2  d\xi
\label{eq:HX5y}\end{equation}
as $N\to\infty$, by the definition of the integral.
Let us here set $N=t^\sigma$  and  compare \e{eq:HX5y}  with \e{eq:HX5z}. First, we obtain relation  \e{eq:HXG}  since the powers of $t$ in 
the left-hand sides of \e{eq:HX5z} and \e{eq:HX5y}  should be the same. It follows that $ \xi^{-r-s/2} $ in \e{eq:HX5x} can be replaced by $ \xi^{-1/2} $. Second, comparing the right-hand sides, using definition \e{eq:HX5x}   and taking into account that
$\| f \|=  \| g\|_{L^2 (\Lambda)}$, we see that
 \begin{equation} 
2\pi \int_{\lambda^{-1}(\Lambda)}  \xi^{-1}    |\omega' (\lambda (\xi) )|^{-1}
\varkappa^2   (\lambda (\xi)) |g(\lambda (\xi)) |^2 d\xi= \int_{\Lambda}    |g ( \lambda)|^2  d\lambda.
\label{eq:HX5Y}\end{equation}
Differentiating relation  \e{eq:HX4A}, we find that
  $
   \omega' (\lambda(\xi))  \lambda'(\xi)=  -s \xi^{-s-1}= -s \xi^{-1}    \omega (\lambda(\xi)) $. 
   Substituting this expression for   $
   \omega' (\lambda(\xi))$ into the left-hand side of \e{eq:HX5Y}, we rewrite  equality \e{eq:HX5Y} as
    \[
2\pi s^{-1}\int_{\Lambda}  \varkappa^2   (\lambda )   \omega (\lambda) ^{-1} |g(\lambda  ) |^2 d\lambda= \int_{\Lambda}    |g ( \lambda)|^2  d\lambda.
\]
Since  $g\in C_{0}^\infty (\Lambda)$ is arbitrary, this yields  relation \e{eq:HXG1} between asymptotic coefficients in \e{eq:HX1}, \e{eq:HX2}  and the spectral measure.

 According to formula \e{eq:HX5}  the functions $(e^{-iJ  t} f )_{n}$ ``live" in the region where $n\sim |t| ^\sigma$. For example, for the Laguerre polynomials,  it follows from \e{eq:Lag4}  that 
 $\sigma=2$. This is fairly unusual since scattering states are normally  concentrated  in the region where $n\sim |t|$. Similarly,   for continuous operators of the Schr\"odinger type we have the relation $x\sim  |t|$.  Indeed, consider, for example, the operator $H=-d^2/ dx^2$ in the space $L^2 ({\Bbb R})$. In this case, we have
                \[
(e^{-i H t} f) (x) =  e^{\mp \pi i/4} (2 |t|)^{-1/2} e^{i x^2/ (4t)} \hat{f} (x/(2t)) + \varepsilon (x,t)
\]
where   $\hat{f}(\xi) =(2\pi)^{-1/2}   
 \int_{-\infty}^\infty e^{-   i x \xi} f( x) d x$ is the Fourier transform of $f(x)$ and the   norm in $L^2 ({\Bbb R})$ of  the term $\varepsilon (\cdot ,t)$  tends to zero as $|t|\to\infty$.
  
  The above arguments relied on the stationary phase method and  strongly used the assumption $\omega' (\lambda)\neq 0$. Let us  now consider, at a very heuristic level, the dispersionless case $\omega' (\lambda) = 0$.  We suppose that condition \e{eq:HX1} is satisfied with $\Omega_{n}  (\lambda)=\nu_n \lambda  + \d_{n}$ where $\lambda\in\Bbb R$, $\nu_{n}  =\omega    n^s + o(n^s)$, $s\in (0,1) $ and $\d_{n}$  do not depend on $\lambda$. Then it  follows from \e{eq:HX3} where $\Lambda={\Bbb R} $ that, up to a  term which tends to zero in $\ell^2 ({\Bbb Z}_{+})$ as $|t| \to \infty$,
   \begin{equation} 
(e^{-i J t} f)_{n} = (2\pi)^{1/2} (n+1)^{-r}   \big( e^{i\d_{n}}\widehat{G}(t-\nu_n) + e^{-i\d_{n}} \widehat{G}(t+\nu_n ) \big) 
\label{eq:HY}\end{equation}
  where   $  \widehat{G}(t)  $
 is the Fourier transform of $G(\lambda)=\varkappa   (\lambda ) g(\lambda)$.  If $\omega>0$ and $t\to+\infty$,  the second term in the right-hand side   is negligible. The operators $e^{-i J t}$  being unitary, it follows from \e{eq:HY}  that
 \begin{equation} 
2\pi  \sum_{n=0}^\infty   (n+1)^{-2r}|\widehat{G}(t-\nu_n )|^2 \to \| f\|^2 = \int_{-\infty}^\infty    |g ( \lambda)|^2  d\lambda \q \mbox{as} \q t \to + \infty. 
\label{eq:HY2}\end{equation}

It is natural to expect that the limit of the left-hand side here is determined by $n$ such that $\omega n^s \sim t$  whence $n^{-2r} \sim (t/\omega)^{-2\sigma r}$.  Let us   set $m= n- (t/\omega)^\sigma$  and use that
$t-\nu_n \sim - s t (\omega /t)^\sigma m$.
Then \e{eq:HY2}  implies that
 \begin{equation} 
\Big(2\pi t^{-2\sigma r+\sigma -1}  \omega^{ 2\sigma r-\sigma} s^{-1}  \Big) \Big( N^{-1}\sum_{m=-\infty}^\infty   |\widehat{G}( -   m /N )|^2 \Big)\to \int_{-\infty}^\infty    |g ( \lambda)|^2  d\lambda
\label{eq:HY3}\end{equation}
where $ N= s^{-1} t^{\sigma-1}\omega^{-\sigma} \to\infty$. According to \e{eq:HX5y}  and the Parseval identity the second factor in the left-hand side has a finite limit   $\int_{-\infty}^\infty     |G ( \lambda)|^2  d\lambda$.  Therefore the power of $t$ in the first factor should be zero which yields equality \e{eq:HXG}.  Now relation \e{eq:HY3} shows that
 \[
2\pi  \int_{-\infty}^\infty     |G ( \lambda)|^2  d\lambda= \omega s \int_{-\infty}^\infty     |g ( \lambda)|^2  d\lambda.
\]
Since  $G(\lambda)=\varkappa   (\lambda ) g(\lambda)$ and $g\in C_{0}^\infty ({\Bbb R})$ is arbitrary, we again arrive at equality \e{eq:HXG1}  where $\omega(\lambda)=\omega$  does not depend on $\lambda$.

Let us summarize the results obtained. Suppose that the orthonormal polynomials  $P_{n}  (\lambda)$ have asymptotic  behavior \e{eq:HX1} with  the phase $\Omega_{n} (\lambda)$ satisfying \e{eq:HX2}.  Then, necessarily relations \e{eq:HXG}  and \e{eq:HXG1}  hold true.  Precise conditions guaranteeing \e{eq:HXG}  and \e{eq:HXG1} and proofs of these relations will be published elsewhere.

        \subsection{Laguerre polynomials} 
        
        Let the operators $J_{p}$ be defined by formula \e{eq:RR} with the coefficients \e{eq:Lag}.
            Theorem~\ref{WW}  shows that, for all $p,q> -1$, 
       \[
\lim_{t\to\pm\infty}\|  e^{-iJ_{p}t} f - e^{-iJ_{q}t}\ti{f}_{\pm} \| = 0 \q \mbox{if}  \q
\ti{f}_{\pm} = W_{\pm} (J_{q} , J_{p})f,
\]
that is,  the time        evolution  of $e^{-iJ_{p}t} f $ is the same for all $p>-1$; only the initial data are changed.

Our goal here is to obtain detailed asymptotic formulas  for $( e^{-iJ_{p}t} f )_{n}$  as $t\to \pm \infty$.
We  choose  $ f$ from the set $\Phi_{p}^* C_{0}^\infty ({\Bbb R}_{+})$  dense in $L^2 ({\Bbb R}_{+})$. It  turns out that the asymptotics of $( e^{-iJ_{p}t} f )_{n}$  depends crucially on the ratio $n/t^2$.
Below we omit the index $p$.  Since $e^{-iJ t} f  
=    \Phi^* e^{-i{\sf A} t}\Phi f$,  Lemma~\ref{WW3}  can be reformulated as follows.
 
 \begin{lem}\label{WU1-}
     Let $g= \Phi f \in C_{0}^\infty ({\Bbb R}_{+})$, and let the operators $V_{\pm} $ be defined by formula \e{eq:WW4}.
     Then 
     \[
e^{-iJ t} f = \nu^{\pm 1}  V_{\pm} e^{-i {\sf A} t} g + \varepsilon_{\pm} (t) 
\]
where   $\| \varepsilon_{\pm} (t)\| \to 0$ as $t\to \pm \infty$. 
     \end{lem}

Thus,  we have to find the asymptotics of 
 \begin{equation}
( V_{\pm} e^{-i {\sf A} t} g)_{n}  = (-1)^n 2^{-1} \pi^{-1/2}  (n+1)^{-1/4}
 \int_{0}^\infty e^{\pm  2 i \sqrt{n\lambda} -i\lambda t} G(\lambda) d\lambda,  
   \label{eq:WW4a}\end{equation}
where $G(\lambda)= \lambda^{-1/4} g(\lambda)$    as $t\to\pm \infty$. 


   The first assertion  shows that these functions     are small as $|t|\to\infty$ both for   relatively  ``small" and ``large" $n$.
   

    \begin{lem}\label{WU1}
     Let the operators $V_{\pm} $ be defined by formula \e{eq:WW4}, and let $G \in C_{0}^\infty ({\Bbb R}_{+})$. Suppose that $\supp G \subset [\lambda_{1}, \lambda_{2}]$ and choose  $\mu_{1}< \lambda_{1}  $, $\mu_{2}>  \lambda_2 $.
Then for all $n\leq \mu_{1} t^2$  
 and $  n\geq \mu_{2} t^2$,  all $k\in {\Bbb Z}_{+}$ and some constants $C_{k}$, we have estimates
   \begin{equation}
\big|  ( V_{\pm} e^{-i {\sf A} t} g)_{n}  \big|\leq C_{k} (\sqrt{n}+ |t|)^{-k},\q \forall t\in{\Bbb R}.
\label{eq:WU2}\end{equation}
     \end{lem}
     
      \begin{pf}
      We proceed from formula \e{eq:WW7A} and
      estimate its right-hand side.
      If $n\leq \mu_{1} t^2$, then
       \[
       | \sqrt{n}\lambda^{-1/2}-t |\geq |t | - \sqrt{n}\lambda_{1}^{-1/2}
       \geq |t | (1- (\mu_{1}/\lambda_{1})^{1/2})\geq c (\sqrt{n}+ |t|) 
       \]
       where $c (\sqrt{\mu_{1}}+1)= 1-\sqrt{\mu_1/\lambda_{1}} > 0$.  Quite similarly, if
       $n\geq \mu_2 t^2$, then
       \[
       | \sqrt{n}\lambda^{-1/2}-t |\geq   \sqrt{n}\lambda_2^{-1/2} - |t | 
       \geq \sqrt{n} (\lambda_2^{-1/2}- \mu_2^{-1/2})\geq c (\sqrt{n}+ |t|) 
       \]
       where $c (\sqrt{\mu_2}+1)=  \sqrt{\mu_2/\lambda_2} - 1 > 0$. According to \e{eq:WW7A} this proves \e{eq:WU2} for $k=1$. 
       
        Integrating by parts  $k$ times in \e{eq:WW7A}, we obtain \e{eq:WU2}.
        \end{pf}
     
   
To find asymptotics of the   integral in \e{eq:WW4a},     we use the stationary phase method.
   Put  
           \[
             \xi= \sqrt{n} t^{-1} \q  \mbox{and}\q
  \theta ( \lambda, \xi)  =  -2\xi \sqrt{\lambda} +  \lambda .
\]
Then
    \begin{equation}
 \int_{0}^\infty e^{   2 i \sqrt{n\lambda} -i\lambda t} G(\lambda) d\lambda
 =  \int_{0}^\infty e^{ -   i  \theta ( \lambda, \xi) t}G(\lambda) d\lambda=: {\cal I}(t, \xi).
   \label{eq:WU4}\end{equation}
Differentiating   $ \theta(\lambda,\xi) $ in $\lambda$, we find that
   \[
 \theta '(\lambda,\xi)  =  - \xi \lambda^{-1/2} +   1\q \mbox{and} \q \theta ''(\lambda, \xi)  = 2^{-1} \xi \lambda^{-3/2} .
  \]
  The stationary point $\lambda_{0} =\lambda_{0}(\xi)$ of the phase  $ \theta(\lambda,\xi)$ is determined by the equation
  $ \theta'(\lambda_{0},\xi)  = 0$  whence $
\lambda_{0} = \xi^2  $.
Therefore the stationary phase method yields
 \[
 {\cal I} (t,\xi) = e^{\mp i \pi/4 } e^{ - i \theta ( \lambda_{0}, \xi)t }\sqrt{\frac{2\pi}{| \theta''(\lambda_{0}, \xi)t |}} \; G(\lambda_{0}) + o(|t|^{-1/2}) 
    \]
   as $t\to \pm\infty$. Since $\theta (\lambda_{0}, \xi)= - \xi^{2}$ and $\theta''(\lambda_{0}, \xi)= 2^{-1}\xi^{-2}$, we arrive at the following intermediary result.

   \begin{lem}\label{WUx}
   Let  $G\in C_{0}^\infty ({\Bbb R}_{+})$.
 Then the integral \e{eq:WU4} has asymptotics
   \begin{equation}
 {\cal I} (t,\xi) =  2 \pi^{1/2}e^{\mp i \pi/4 } e^{i   \xi^2 t }  | t |^{-1} |\xi | G(\xi^2)  + o(|t|^{-1}) , \q t\to \pm\infty,
   \label{eq:WU5+}\end{equation}
   with the estimate of the remainder uniform in $\xi$ from  compact subintervals of ${\Bbb R}\setminus\{0\}$.
 \end{lem}
 
  Let us come back to formula \e{eq:WW4a}. Set
       \begin{equation}
(U(t) g )_{n}=   (-1)^n e^{in/t} |t|^{-1} g (n/t^2), \q t\neq 0.  
\label{eq:KK}\end{equation}

\begin{lem}\label{WUy}
 Let  $G\in C_{0}^\infty ({\Bbb R}_{+})$ and $0<\mu_{1}  <\mu_{2}  <\infty$.
Then 
 \begin{equation}
 \sup_{n\in (\mu_{1}t^2, \mu_2 t^2)}\big| (V_{\pm}e^{-i {\sf A} t} g)_{n}-e^{\mp \pi i/4}(U(t) g )_{n}\big| = o(|t|^{-1}), \q t\to\pm \infty.
   \label{eq:KK1}\end{equation}
 \end{lem}
 
  \begin{pf}
  Let us set in \e{eq:WU5+} $G(\lambda)= \lambda^{-1/4}  g(\lambda)$ and  $\xi=\sqrt{n}t^{-1}$ so that
  \[
  |\xi |G(\xi^2)= |\xi |^{1/2} g(\xi^2)=  n^{1/4} |t|^{-1/2} g(n/t^2).
  \]
  Thus,  it follows from Lemma~\ref{WUx}  that
 \[
  \int_{0}^\infty e^{  \pm 2 i \sqrt{n\lambda} -i\lambda t} \lambda^{-1/4} g (\lambda) d\lambda
  = 2 \pi^{1/2}e^{\mp i \pi/4 } n^{1/4} e^{i   n/t  }    |t|^{-1}  g(\frac{n}{t^2}) + o(|t|^{-1}),\q t\to\pm\infty,
  \]
  as long as $n\in (\mu_{1}t^2, \mu_2 t^2)$.
  Putting together this relation with \e{eq:WW4a}, we arrive at \e{eq:KK1}.
    \end{pf}

  Now we are in a position to obtain an asymptotic formula for $e^{-iJ_{p} t} f$ as $t\to\pm \infty$. 
   
   \begin{thm}\label{WU}
   Let $J_{p}$, $p> -1$,  be  a Jacobi operator with matrix elements \e{eq:Lag}, and let  the operator
    $U (t)$ be given by formula \e{eq:KK}. Define the operators $\Phi_{p}$   by formulas \e{eq:UF} and \e{eq:UU1} and  
   suppose that   $\Phi_{p} f\in C_{0} ^\infty ({\Bbb R}_{+} )$. Then 
     \begin{equation}
\lim_{t\to\pm \infty} \| e^{-iJ_{p} t} f - e^{\mp i (p+1)\pi/2} U (t)\Phi_{p} f\| =0.
\label{eq:WU}\end{equation}
 \end{thm}
 
   \begin{pf}
 According to Lemma~\ref{WU1-} we  can replace here $e^{-iJ t} f $  by $ \nu^{\pm 1}  V_{\pm} e^{-i {\sf A} t} g$  where 
      $g= \Phi f  $. Suppose that $\supp g \subset [\lambda_{1}, \lambda_{2}]$ and choose
      $\mu_{1}< \lambda_{1}  $,  $\mu_{2}>  \lambda_2 $.  It follows from Lemma~\ref{WU1} that
      \begin{equation}
      \big(\sum_{n\leq \mu_{1}t^2}+ \sum_{n\geq \mu_2 t^2}\big) |(V_{\pm}   e^{-i {\sf A} t} g)_{n} |^2\to 0
\label{eq:KK2}\end{equation}
      as $t\to\pm \infty$.  According to  \e{eq:KK1} we also have
      \[
     \sum_{\mu_1 t^2\leq n\leq \mu_{2}t^2}|( V_{\pm}e^{-i {\sf A} t} g)_{n}-e^{\mp \pi i/4}(U(t) g )_{n}|^2 = o(1 ) .
      \]
      Combined with \e{eq:KK2} this yields relation \e{eq:WU}.
   \end{pf}


According to \e{eq:Lag4} for the Laguerre polynomials, we have $\Lambda={\Bbb R}_{+}$, $r=1/4$, $s=1/2$ and
\[
\kappa(\lambda)=  \frac{1} {2}  \sqrt{\frac{\Gamma (1+p)  }
 {\pi  }}  \lambda^{-p/2-1/4}  e^{\lambda/2} , \q
\Omega_{n} (\lambda)=\pi n+ 2\sqrt{n\lambda}-\frac{2p+1} {4} \pi, \q \omega(\lambda)=\lambda^{-1/2}.
 \] 
Since $\tau(\lambda)$ is given by  \e{eq:Lag2},  the identity \e{eq:HXG1} is satisfied.


            \subsection{Jacobi polynomials} 
            
            In this subsection we define a Jacobi operator   by its spectral  measure $d\rho (\lambda)$. We suppose that this measure is supported on the interval $[-1,1]$ and  
   \begin{equation}
d\rho (\lambda)= \tau (\lambda) d\lambda, \q  \lambda\in (-1,1),
\label{eq:Jac}\end{equation}
where
 \begin{equation}
 \tau (\lambda)= k (1-\lambda)^{\alpha} (1+\lambda)^{\beta},  \q \alpha ,\beta>-1.
\label{eq:Jac1}\end{equation}
The weight function $ \tau(\lambda)= \tau_{\alpha ,\beta}(\lambda)$ (as well as all other objects discussed below) depends on $\alpha $ and $\beta$, but these parameters are often omitted in notation. 
 The constant 
$  k = k_{\alpha ,\beta}$
 is chosen in such a way that the measure \e{eq:Jac} is normalized, i.e., 
$\rho ({\Bbb R})=\rho ((-1,1))=1$. 
The orthonormal polynomials ${\sf G}_{n}(\lambda)= {\sf G}_{n}^{(\alpha,\beta)}(\lambda)$ determined by the measure
\e{eq:Jac}, \e{eq:Jac1} are known as the Jacobi polynomials.

Let $J=J_{\alpha,\beta}$ be the Jacobi operator  with the spectral measure $ d\rho(\lambda)= d\rho_{\alpha ,\beta}(\lambda)$. Explicit expressions for
its matrix elements  $a_{n}, b_{n}$ can be found, for example, in the books \cite{BE, Sz}, but we do not need them. We here note only  asymptotic  formulas 
 \begin{equation}
 a_{n}  =1/2 + 2^{-4} (1 -2\alpha^2-2\beta^2) n^{-2} +O\big( n^{-3}\big),\q
b_{n}  = 2^{-2}  (\beta^2-\alpha^2) n^{-2} +O\big( n^{-3}\big) 
\label{eq:norm6}\end{equation}
for  the matrix elements and  (see formula (8.21.10)   in the book~\cite{Sz})
  \begin{multline}
 {\sf G}_{n} (\lambda)=  2^{1/2}( \pi k)^{-1/2}(1-\lambda)^{-(1+2\alpha)/4}(1+\lambda)^{-(1+2\beta)/4} 
 \\
 \times \cos \big((n+ \gamma) \arcsin\lambda - \pi(2n+\beta-\alpha) /4\big) + O (n^{-1}), \q \gamma = (\alpha+\beta+1)/2,
\label{eq:GG} \end{multline}
 for the orthonormal polynomials.  Estimate of the remainder in \e{eq:GG} is uniform in $\lambda$ from compact subsets of $(-1,1)$.

   Similarly to the cases of the Laguerre   polynomials, we set
  $
\varphi_{n} (\lambda) =   \sqrt{\tau(\lambda)}  {\sf G}_{n} (\lambda)
$
   and define the mapping $\Phi: \ell^2 ({\Bbb Z}_{+})\to L^2 (-1,1)$ by   formula \e{eq:UU1} where $\lambda\in (-1,1)$.
   Using \e{eq:Jac1} and  \e{eq:GG}, we obtain  the representation 
 \begin{equation}
   e^{-i J t} f= \Phi^* e^{-i{\sf A}t}g =      V_{+} e^{-i{\sf A}t}g
    +   V_{-} e^{-i{\sf A}t}g  + R  (t) g, \q g= \Phi f \in C_{0}^\infty ( -1, 1),
   \label{eq:H6}\end{equation}
where
   \[
( V_{\pm} g)_{n}  =     (2\pi)^{-1/2} i^{\mp n} e^{\pm i ( \alpha -\beta)\pi/4} 
 \int_{-1}^1 e^{\pm   i  (n+ \gamma) \arcsin \lambda } (1-\lambda^2)^{-1/4}g(\lambda) d\lambda,
 \]
  and the remainder $R  (t) g$ satisfies  condition \e{eq:WW61}.  

Let us state analogues of Lemmas~\ref{WW2} and \ref{WU1}.

\begin{lem}\label{WW2J}
     For an arbitrary $G\in C_{0}^\infty  ( -1, 1)$,  an estimate
   \begin{equation}
\Big|  \int_{-1}^1 e^{\pm   i  n \arcsin \lambda -i\lambda t} G(\lambda) d\lambda \Big|\leq C_{k} (n+ | t| )^{-k} \label{eq:WW7J}\end{equation}
is true for all $k\in {\Bbb Z}_{+}$ if $\mp t>0$.  
    \end{lem}
    
     \begin{pf}
     Integrating by parts, we see that
   \[
     \int_{-1}^1 e^{\pm   i  n \arcsin \lambda -i\lambda t}G(\lambda) d\lambda 
          = 
 \pm  i       \int_{-1}^1 e^{\pm   i  n \arcsin \lambda -i\lambda t}\Big(\frac{G(\lambda)}{ n (1-\lambda^2)^{-1/2}  \mp t}\Big)' d\lambda  .
 \]
    This yields \e{eq:WW7J} for $k=1$ because $ \mp t=|t|$ and $|\lambda|\leq c<1$ on the support of $G$.
      Further    integrations  by parts     lead to  estimates \e{eq:WW7J} for all $k$.
        \end{pf} 
        
       Quite similarly, we obtain  also the following result.
     
\begin{lem}\label{WU1J}
   Let $G\in C_{0}^\infty  (( -1, 1)\setminus\{0\})$.
Then estimates  \e{eq:WW7J} hold for all $t\in{\Bbb R}$ and all $k\in {\Bbb Z}_{+}$ if either $n\leq \d |t|$  or $  n\geq (1-\d) | t |$  for    a sufficiently small   number $\d$  depending on $\supp G$.
     \end{lem}

     We use these lemmas with 
           \begin{equation}
            G_{\pm}(\lambda)=e^{\pm   i    \gamma \arcsin \lambda } (1-\lambda^2)^{-1/4}g(\lambda) 
             \label{eq:GG2y}\end{equation}
             and take equality \e{eq:H6}  into account.
             According to Lemma~\ref{WW2J} relation \e{eq:WW7+} is satisfied so that it   suffices to consider $V_{\pm}
     e^{-i{\sf A}t}g$ as $t\to\pm\infty$. According to Lemma~\ref{WU1J} we only have to study
      \[
( V_{\pm}    e^{-i{\sf A}t} g)_{n}  =     (2\pi)^{-1/2} i^{\mp n} e^{\pm i ( \alpha -\beta)\pi/4} 
 \int_{-1}^1 e^{\pm   i  n \arcsin \lambda-i\lambda t } G_{\pm}(\lambda) d\lambda
   \]
 for $\d |t| \leq  n\leq (1-\d) | t |$  where $\d>0$.

Let us now set
  \begin{equation}
  \xi= n |t|^{-1},\q
 \theta ( \lambda, \xi)  =   \xi \arcsin\lambda  -  \lambda 
\label{eq:GG4}\end{equation}
and consider an integral 
  \begin{equation}
  {\cal I}(t, \xi)
 =  \int_{-1}^1 e^{    i   \theta ( \lambda, \xi) t}G (\lambda) d\lambda , \q \xi\in (0, 1),
   \label{eq:GG3}\end{equation}
   where  $G \in C_{0}^\infty (( -1, 1)\setminus\{0\}) $.
Differentiating   $ \theta (\lambda,\xi) $ in $\lambda$, we see that
   \[
\theta '(\lambda,\xi)  =   \xi (1- \lambda^2)^{-1/2}  -  1\q \mbox{and} \q \theta ''(\lambda, \xi)  =   \xi \lambda (1-\lambda^2)^{-3/2} .
  \]
  The stationary points $\lambda_\pm  =\lambda_\pm(\xi)$ of the phase  $ \theta(\lambda,\xi)$ are determined by the equation
  $\theta '(\lambda_\pm ,\xi)  = 0$  whence
       \[
       \lambda_\pm=\pm\lambda_0, \q
\lambda_0  =   \sqrt{1-\xi^2 }  
\]
and  
\[
\theta''(\lambda_{\pm}, \xi)= \pm    \xi^ {-2}\sqrt{1-\xi^2}  .
\]
Note also  that
 \begin{equation}
 \theta (  \lambda_{0} , \xi)= \xi \arccos\xi  - \sqrt{1-\xi^2}=: \psi  (\xi) .
\label{eq:WNJ}\end{equation}

Therefore the stationary phase method yields the following  intermediary result.

   \begin{lem}\label{GG}
   Let the phases $ \theta ( \lambda, \xi)$  and $\psi  (\xi)$  be  given  by formulas   \e{eq:GG4}
   and \e{eq:WNJ}, respectively.
 Then the integral \e{eq:GG3}  has asymptotics 
 \begin{multline}
 {\cal I} (t,\xi) =  (2\pi)^{1/2}|t|^{-1/2} \xi (1-\xi^2)^{-1/4} 
 \\
 \times  \Big( e^{\pm i \pi/4 } e^{ i \psi  (\xi)t} G \big(\sqrt{1-\xi^2}\,\big) + e^{ \mp i \pi/4 } e^{ - i \psi  (\xi) t} G \big(- \sqrt{1-\xi^2}\,\big) \Big)+ o(|t|^{-1/2}) 
    \label{eq:GG6}\end{multline}
       as $t\to \pm \infty$.
  The estimate of the remainder in \e{eq:GG6} is uniform in $\xi$ from  compact subintervals of $  (0,1) $.
 \end{lem}
 
 Let  us  now  apply Lemma~\ref{GG} to functions  \e{eq:GG2y}. We set
 \begin{equation}
 h_{\pm}  (\xi)=e^{\pm i \pi/4 } e^{\pm   i    \gamma \arccos \xi  }    \xi^{1/2}  (1-\xi^2)^{-1/4} g\big( \pm \sqrt{1-\xi^2} \big)
 \label{eq:KKJx}\end{equation}
 and
 \begin{equation}
(U(t) g )_{n}=  i^{\mp n}  e^{\pm i ( \alpha - \beta)\pi/4}    |t|^{-1/2}   
 \big( e^{ i \psi  (\xi)t}  h_{+}  (\xi)+ e^{- i \psi  (\xi)t}  h_{-}  (\xi)\big)    
\label{eq:KKJ}\end{equation}
where $ \xi= n |t |^{-1} \in (0,1)$, $\pm t > 0$ and  the function $\psi 
(\xi) $ is defined by formula \e{eq:WNJ}.  For $n\geq  |t|$, we set $(U(t) g )_{n}=0$.

      Putting the results obtained together, we state our final result.

 \begin{thm}\label{WJJ}
Let $J=J_{\alpha,\beta}$ be the Jacobi operator  corresponding to the weight function  \e{eq:Jac}. Define the operator
      $\Phi: \ell^2 ({\Bbb Z}_{+})\to L^2 (-1,1)$ by   formula \e{eq:UU1}     and  
   suppose that   $g=\Phi f\in  C_{0}^\infty   (( -1, 1)\setminus\{0\})$. Let  the operator
    $U (t)$ be given by equalities \e{eq:KKJx}  and \e{eq:KKJ}.
 Then 
     \[
\lim_{t\to\pm \infty} \| e^{-iJ  t} f -  U (t)\Phi  f\| =0.
\]
 \end{thm}
 
 The operator $J_{1/2,1/2}=: J^{(0)}$ is particularly simple. It is known as the free discrete Schr\"odinger operator. In this case,  we have $a_{n}=1/2$, $b_{n}=0$ for all $n$.  The corresponding orthonormal polynomials ${\sf G}_{n}^{(1/2 , 1/2)}(\lambda)$ are the Chebyshev polynomials of second kind.  According to \e{eq:norm6}   the difference $J- J^{(0)}$ is trace class for all $\alpha$ and $\beta$, so that the wave operators $W_{\pm} (J, J^{(0)})$ as well as $W_{\pm} ( J^{(0)},J)$ exist. Similarly to Theorem~\ref{WW}, this fact can be also checked by a direct calculation using relations  \e{eq:H6} and \e{eq:WW7+}; see \cite{Y/LD} for details.   This yields also explicit expressions for the wave operators:
 $
 W_{\pm} (J, J^{(0)})=\Phi^*\Sigma_{\pm} \Phi^{(0)}
$
 where $\Sigma_{\pm}$ is the multiplication operator in the space $L^2 (-1,1)$  by  the function
\[
  e^{\mp i \pi(  \alpha-\beta)/4}
  e^{\mp i 2^{-1}  (\alpha+\beta-1) \arcsin\lambda}  .
   \]
 
 For Jacobi polynomials, we have $\Lambda= (-1,0)$ or $\Lambda= (0,1)$, $r=0$, $s=1$ and $\omega(\lambda)= (1-\lambda^2)^{-1/2}$. Comparing \e{eq:Jac1} and \e{eq:GG}, we also see that $\varkappa(\lambda)= (2\pi)^{-1/2} (1-\lambda^2)^{-1/4}$ whence $2\pi  \varkappa^2(\lambda)= \omega(\lambda)$ which is consistent with  \e{eq:HXG1}.
 
 We finally note that the asymptotic behavior of the evolution $e^{-iJ  t} f $ looks similar to the corresponding result
 for the Klein-Gordon  operator $\sqrt{-d^2/dx^2 +1}$ in the space $L^2 ({\Bbb R})$.
 
  \subsection{Hermite polynomials} 
      
      The Hermite polynomials $H_{n}  (z)$ are determined by the recurrence coefficients \e{eq:H}.
      As usual,    relations  \e{eq:Py} for $H_{n}  (z)$ are complemented by the   boundary conditions $H_{-1 } (z) =0$, $H_{0 } (z) =1$. 
      According to asymptotic formula (10.15.18) in \cite{BE}, we have
\begin{equation}
H_{n}  (\lambda)=
   2^{1/2}
  \pi^{-1/4} e^{\lambda^2/2}  (2n+1)^{-1/4} \cos \big( \sqrt{2n+1} \, \lambda-  \pi  n/2\big) +O(n^{-3/4})  
\label{eq:H2}\end{equation}
as $n\to\infty$.    This  asymptotics   is uniform in $\lambda\in \Bbb R$   on compact subintervals.
      
  Let us  consider the Jacobi operators $ J $ defined by formula \e{eq:RR} where $a_{n}$, $b_{n}$ are  given by \e{eq:H}.       The spectral measure of $J$ equals  
   $ d\rho  (\lambda)=   \pi^{-1/2} e^{-\lambda^2} d\lambda$ where $ \lambda\in {\Bbb R}$
 (see, e.g., formula (10.13.1) in \cite{BE}).
Thus, $d\rho  (\lambda)$ is absolutely continuous and its 
    support is the whole axis  $\Bbb R$.

  Following the scheme exposed in Sect.~2.2, we reduce the Jacobi operator $J $ to the operator  ${\sf A}$ of multiplication by $\lambda$ in the space $L^2 ({\Bbb R})$.  To that end, we 
   introduce a mapping $\Phi  : \ell ^2 ({\Bbb Z}_{+})\to L^2 ({\Bbb R}  )$  by the formula  
  \begin{equation}
 (\Phi f)(\lambda)= \pi^{-1/4}\sum_{n=0}^\infty e^{-\lambda^2/2} H_{n}  (\lambda) f_{n}, \q f=\{ f_{n}\}_{n=0}^\infty \in {\cal D}, \q \lambda\in {\Bbb R} .
\label{eq:H4}\end{equation}
The operator $\Phi $ is unitary and enjoys  the intertwining property $ \Phi  \, J  = {\sf A} \, \Phi   $.

Putting together formulas \e{eq:H2} and \e{eq:H4}, we obtain representation \e{eq:H6}
   where
     \[
( V_{\pm} g)_{n}  = i^{\mp n}    (2\pi)^{-1/2}  (2n+1)^{-1/4}
 \int_{-\infty}^\infty e^{\pm   i \sqrt{2n+1} \,\lambda } g(\lambda) d\lambda
   \]
            and the remainder $R  (t) g$ satisfies condition \e{eq:WW61}.  Then
   \[
( V_{\pm}  e^{-i{\sf A}t}g)_{n}  =   i^{\mp n} (2n+1)^{-1/4}    \hat{g}(t \mp \sqrt{2n+1})
    \]
    where $\hat{g}(x)$ is  the Fourier transform of $g$.
    Since $  \hat{g}(x) = O( x|^{-k})$ as $|x| \to \infty$ for all $k>0$,  we see that
    \[
\lim_{t\to \pm\infty}    \sum_{n=0}^\infty   (2n+1)^{-1/2}|   \hat{g}(t \pm \sqrt{2n+1})|^2 =0, 
\]
whence relation \e{eq:WW7+} follows.  Thus, representation \e{eq:H6}  implies the result below.

     \begin{thm}\label{WH}
   Let $J $  be  a Jacobi operator with matrix elements \e{eq:H}, and let  the operator
    $U (t)$ be given by the formula 
      \begin{equation}
( U(t) g)_{n}  =   i^{\mp n} (2 n+1)^{-1/4}    \hat{g}(t  \mp \sqrt{2n+1}), \q \pm t>0.
\label{eq:WH}\end{equation}
    Define the operator $\Phi $   by formula \e{eq:H4} and  
   suppose that   $\Phi  f\in C_{0} ^\infty ({\Bbb R} )$. Then 
     \[
\lim_{t\to\pm \infty} \| e^{-iJ  t} f -   U (t)\Phi  f\| =0.
\]
 \end{thm}

     \begin{cor} 
     For all $g \in  C_{0}^\infty (\Bbb R)$, we have
\begin{equation}
\lim_{t\to \pm \infty}    \sum_{n=0}^\infty   (2 n+1)^{-1/2}|  \hat{g}( t\mp \sqrt{2n+1} )|^2 = \| \hat{g}\|_{L^2 ({\Bbb R})}^2  . 
\label{eq:WH1}\end{equation}
\end{cor}

According to \e{eq:WH} the evolution $e^{-iJ t} f$ is   dispersionless. Clearly, it is similar to time evolutions for first order  differential operators.

Finally, we note that for the Hermite polynomials, $\Lambda=\Bbb R$, $r=1/4$, $s=1/2$, $\Omega_{n} (\lambda)=\sqrt{2n+1}\, \lambda -\pi n/2$, $\omega(\lambda)=\sqrt{2}$, $\kappa(\lambda)=  2^{-3/4}\pi^{-1/4} e^{\lambda^2/2}$. Thus, the identity \e{eq:HXG1}  remains true.  Evolution  \e{eq:WH}  is dispersionless now, and relation \e{eq:WH1}  is  a particular case of \e{eq:HY2}.

 \begin{acknowledgments}
Supported by      RFBR grant No.   17-01-00668 A.
\end{acknowledgments}

\small
 


 \end{document}